\documentclass[12pt,a4paper]{article}
\usepackage{cite}
\usepackage{natbib}
\usepackage{breakcites}
\usepackage[normalem]{ulem}
\usepackage{comment}
\usepackage[font=footnotesize,labelfont=bf]{caption}

\let\rho=\varrho
\def\fref#1{Fig.~\ref{#1}}

\def\cref#1{Condition~\ref{#1}}
\def\Cref#1{Corollary~\ref{#1}}
\def\eref#1{(\ref{#1})}

\usepackage{amsmath}
\usepackage{amsfonts}
\usepackage{graphicx}
\usepackage{times}
\usepackage{amsthm}
\usepackage{ amssymb }
\usepackage{color}
\usepackage{mhequ}
\usepackage{dsfont}
\usepackage[scanall]{psfrag}
\usepackage[margin=2.5cm]{geometry}
\usepackage{color}
\usepackage{url}
\usepackage{lipsum}
\usepackage{subfig}
\usepackage{mhequ}

\usepackage{tikz}
\usepackage[graphics, active, tightpage]{}

\def\thecomma{\ifx,\thenewxt \else\ifx;\thenext \else\ifx.\thenext
	\else\ifx!\thenext \else\ifx:\thenext\else\ifx)\thenext \else \
	\fi\fi\fi\fi\fi\fi}
\def\condblank{\futurelet\thenext\thecomma}
\def\ie{{\it i.e.,}\condblank}

\numberwithin{equation}{section}

\theoremstyle{definition} 

\usepackage{stmaryrd}

\let\kappa=\varkappa
\let\phi=\varphi

\def\Rsch{r_{\rm sch}}
\def\d{{\rm d}} 

\def\OO{{\mathcal O}}
\def\integer{{\mathbb Z}}
\def\real{{\mathbb R}}

\def\p2t2{{\tilde p_2^{\,2}}}

\definecolor{bittersweet}{rgb}{1.0, 0.44, 0.37}

\let\epsilon=\varepsilon
\usepackage{authblk}
\begin{document}
\title{The detection of relativistic corrections in cosmological N-body simulations}
\author[1,2]{Jean-Pierre Eckmann}
\author[2]{Farbod Hassani}
\affil[1]{D\'epartement de Physique Th\'eorique and Section de
  Math\'ematiques, University of Geneva, Switzerland}
\affil[2]{D\'epartement de Physique Th\'eorique, University of Geneva, Switzerland}
\maketitle
\abstract{Cosmological N-body simulations are done on massively
  parallel computers. This necessitates the use of simple time
  integrators, and, additionally, of mesh-grid approximations of the
  potentials. Recently,  \cite{Adamek:2015eda, Barrera-Hinojosa:2019mzo} have developed general relativistic N-body simulations to capture relativistic effects mainly for cosmological purposes. We therefore ask whether, with the
  available technology, relativistic effects like perihelion advance can be detected
  numerically to a relevant precision.  We first
study the spurious perihelion shift in the Kepler problem, as a function of the integration method used, and then as a function of an additional interpolation of forces on a 2-dimensional lattice. This is done for
several choices of eccentricities and semi-major axes. Using these results, we can predict
which precisions and lattice constants allow for a detection of the relativistic perihelion
advancein N-body simulation. We find that there are only
  small windows of parameters---such as eccentricity, distance from the central object and
  the Schwarzschild radius---for which the corrections can be detected
  in the numerics.}

\section{Introduction}
We consider here so-called (cosmological) N-body simulations such as
\citep{Adamek2016, Springel2005, Teyssier:2001cp}. In these numerical
studies, potentials between the (many) particles are computed on a
lattice (mesh-grid) because of the way such calculations are implemented on
supercomputers. Additionally, some of these projects (such as
\cite{Adamek:2015eda, Barrera-Hinojosa:2019mzo}) add relativistic
corrections to the forces and therefore to the trajectories of particles. The aim of our study is to give bounds on
the detectability of these effects, given the computational
restrictions of these large-scale projects. We will see that in many
current simulations, the necessary precision to detect relativistic
effects on the orbits of particles can simply not be achieved.

It is of course not difficult to devise codes which will compute the
perihelion advance under relativistic corrections to arbitrary high
precision. It is not the aim of our paper to study such algorithms,
but rather, to see how well the integration algorithms work in the
N-body simulations. In these simulations, because one considers
essentially a gas of many particles, the user is restricted to rather
standard integration methods, which just use the differential
equations, but necessarily can not make use of the many invariants
known for the (non-relativistic) Kepler problem, see {\it{e.g.}}, \citep{Preto:2009cp}. Therefore, we need to
first study the performance of standard integration schemes, such as
Euler, Runge-Kutta, and Leap-Frog, because these are the methods which
are widely used. We will see that only with very high precision one is able to detect
the (usually quite small) relativistic corrections. Once this has been
done, we can turn our attention to the effects of the discretizations
(of space), which give then bounds on the necessary grid constants for
which relativistic effects could be detected. We will determine
parameter regions where the relativistic effects can be detected, and
show that, most often, these regions are quite small.


By starting with the simple Kepler problem, in $\real^2$, we
can concentrate on the different numerical effects in a systematic
and clean way. Even so, the reader should realize that there are
several quantities to be considered. The first is the numerical
precision of the time integrator. We study it here in the context of the
subroutines in ODEX \citep{Hairer1993}, and we also compare it to
other methods, such as Euler,  leapfrog (Verlet-St\"ormer) \citep{Hairer2003} or Runge-Kutta with
fixed time step.

To do this, we quantify the numerical errors on trajectories of
particles revolving around a central object.
This will allow us to give conditions which ascertain
which orbits in a specific
N-body simulation are precise enough to be able to measure the
general relativistic perihelion\footnote{We use ``perihelion'' even if
  the central mass is not the sun.}  shift.

After this, we consider the
particle-mesh N-body scheme, as is widely used,
see {\it e.g.}, \citep{Adamek2016, Springel2005}. In it,  forces (coming from fields and potentials)
are discretized and represented on a lattice. Such
elements \citep{arnold2002} are then used
to compute the values of the fields at the
particles' positions.

The force interpolation approximations are usually piecewise
differentiable, and, depending on the implementations mentioned above,
use different elements. It is clear that if the mesh size of the
approximation (of the force) goes to 0, so will the error. But the
relevant question here is to quantify what kind of phenomena
can be captured, given the numerous hardware and modeling constraints.

Particle-mesh N-body simulations are used to study the evolution of
particles under gravity. These codes can be used to study
systems at different scales, from cosmological scales to the size of the
solar system, as the methods and forces are appropriate for all
scales.
In the particle-mesh N-body scheme
\citep{Springel2005,Adamek2016}, space-discret\-izations are performed
to take care of the large
number of particles.  We analyze two common force
interpolations which are used for N-body simulations purposes, namely
the so-called 
linear and bilinear methods,
which respectively correspond to first and second order interpolation.
We will see
that, under conditions to be 
specified, the effect of spatial discretization can be quite large and sometimes
depends on the angle $\theta$ between the direction of the perihelion
and the axes of the discrete lattice. This happens when
the discretization produces discontinuous forces, \ie for the first
order force interpolation. In this case the maximal errors are proportional to
the lattice constant $dx$.
On the other hand, in the second order force interpolation, when one
varies $\theta$, there is a
small perihelion shift, fluctuating around 0.
These fluctuations are
seen to be of order $dx^{1.3}$. Due to the highly
nonlinear step size of ODEX, we are not able to derive analytically
this size of the fluctuations.

Our numerical tests show that, unless the discretization is 
extremely fine, the system will show an uncertainty of the perihelion,
for the Kepler problem, for both force interpolation methods. Our
calculations give limits on the detectability of relativistic effects,
as a function of method, lattice spacing,  as
well as eccentricity and the relativistic parameter  $\Upsilon \equiv
r_{\rm sch}/ r_{\rm per}$, the ratio of perihelion distance of the orbit
$r_{\rm{per}}$ and the Schwarzschild radius of the central mass $r_{\rm
  sch}.$

\section{Using standard time integrators}

Our main interest is the detectability of general relativistic effects in
N-body simulations, and in particular the study of discretization
effects.
But we first need to be sure that the time integration which is used
in these projects  does not already destroy
the precision of the result more than the effect of the space
discretization.
This is the subject of this section.

A Hamiltonian problem can be integrated either as a motion in
Euclidean space, or one can exploit the underlying symplectic
structure of the problem. 
Of course, there are very good symplectic
integrators, \citep{Hairer2010}, but we decided not to use them, for
two reasons.

The first is that in the particle-mesh N-body codes, the  Euclidean
approach is used to solve the system including particles and the
fields on the lattice coordinates and expressing the system in
symplectic coordinates is difficult.  Second, as noted in
\citep{Hairer2003} even the symplectic methods do not preserve the
Runge-Lenz-Pauli vector (the
orientation of the semi-major axis). This means that
because of our focus on the relativistic perihelion advance, the
symplectic integrators present no particular advantage. 
So we will stick with the 
classical high-order Runge-Kutta integrators ODEX \citep{Hairer1993}.
Because it
allows for ``continuous output'' we can use it to determine easily the
advance of the perihelion with high precision.

We summarize those properties of ODEX which are relevant for our study.
As we
will be working with elements to interpolate forces, we need to
explain how the algorithm deals with discontinuities. This is
illustrated in
\citep[Chapter II.9 and in particular, Fig.~9.6]{Hairer1993}.
In the
interior of a plaquette, the algorithm chooses a high enough order to reach
the required tolerance with a large time step $h$. On approaching the
singularity, the algorithm lowers the order (to 4) but decreases
$h$. In fact, the jump is approximated by a polynomial of degree 4,
and this defines something like a new initial condition across the
discontinuity. In the case of Runge-Kutta with fixed time step, the
paper \citep{Back2005} shows that there is a mean systematic error
across the jump, which can be viewed as a weighted combination of
evaluations of the vector field across the singularity.

\section{Perihelion variation for time integrators with fixed time step}

Here, we study the precision of perihelion calculations when working
with \emph{exact} forces. In later sections, we will then see how the
grid discretization further affects this precision. Certainly, if we
want to discriminate between non-relativistic and relativistic
effects, already the integration with exact forces needs to be precise
enough. This will force us to choose a small enough time step $h$.

We analyze the perihelion shift for several 
standard time integrators with fixed time step, namely Euler,
Newton-St\"ormer-Verlet-leapfrog, 2nd order and 
4th order Runge-Kutta.

We call the time step $h$ and we solve the Kepler
problem  in the form\footnote{All positions, velocities, and the like are in $\real^2$.}
\begin{equ}\label{eq:newton}
  \ddot{ x}(t) =   F(t, x)/m~,
\end{equ}
with $x$, $v$,
and $F$ in $\real^2$ and $m$  the mass of the object.
For the convenience of the reader, we spell out these well-known methods.
\paragraph{Euler method}
In this case, we solve \eref{eq:newton} in the form
 \begin{equa}
   {v}_{ n+1}  &=    v_{ n} +   a_{ n} {h}~,\\
    x_{ \rm n+1} &=   x_{ n}+   v_{ n+1} h~,
\end{equa}
where $ v_{ n}$ is the velocity vector at time step $({ n}$)
which is defined as $ v_{ n} \equiv ({ x_{ n} - 
  x_{ n-1} })/{h}$ and $ a_{ n} $ is the acceleration and is
defined as the ratio of the force and mass $ a_{ n}\equiv  
a(t_{ n},  x_{ n}) =   F (t_{ n},  x_{
  n})/m$. 
It is well-known that this implicit/backward Euler method is more
stable than the explicit/forward Euler method. But it is somewhat more
difficult to implement
for non-linear differential equations.

\paragraph{Newton-St\"ormer-Verlet-leapfrog method }
For this widely used method (sometimes called
 ``kick-drift-kick'' form of leap-frog) \citep{Hairer2003}[Eq.~1.5], the updates are
\begin{equa}
  v_{n+1/2}  &=  v_{n}  + a_{n} \frac{h}{2}~,\\
 x_{n+1} &= x_{n} +  v_{n +\frac 12} h~,\\
    v_{n+1}  &=  v_{n + \frac 12}  +  a_{n+1} \frac{h}{2} ~.
\end{equa}
This method is used more often as it is a symplectic method and stable
and is shown to work very well for various stiff ODEs \citep{Hairer2003}.
\paragraph{Second order Runge-Kutta method}
Finally, we will comment on the perihelion advance for the 2nd
and  4th order Runge-Kutta methods (with fixed time step $h$). The Kepler problem using 2nd order Runge-Kutta algorithm---also known as midpoint method---reads,
\begin{equa}
 k^{(1)}_{ x} &= {v}_{n} ~,\\
 k^{(1)}_{ v} &= {a}_{n} ~,\\
 k^{(2)}_{ x} &= {v}_{n + \frac12} ~,\\
 k^{(2)}_{ v} &= {a}_{n + \frac12} ~,\\
 x _{n+1}&=  x _{n} +  k^{(2)}_{ x} h ~,\\ 
 v _{n+1}&=  v_{n} +  k^{(2)}_{ v} h ~,\\ 
\end{equa}
where $k^{(1)}_{ x}$ is the estimate of velocity (derivative of $ x$) in time step $\rm n$, $ k^{(1)}_{ v}$ is the estimate of acceleration (derivative of $ v$) in time step $\rm n$ and the same for  $ k^{(2)}_{ x}$ and $ k^{(2)}_{ v}$. The acceleration at time $ n+ \frac12$  is obtained by
\begin{equa}
 {a}_{n + \frac12} \equiv  \frac{F (t_{n +\frac12},  x_{n+\frac12})}{m}  = \frac{ F \big(t_{n +\frac12},  x_{n}+  k^{(1)}_{ x}  h/2 \big)}{m}  ~.
\end{equa}
Also, to obtain the velocity at time $ n+ \frac12$ we need to use $ k^{(1)}_{ v} $. The corresponding tableau for the second order Runge-Kutta method for each first order differential equation is
\begin{equa}
\begin{tabular}{c|c c}
   0 &  &  \\
   1/2 & 1/2 & \\
         \hline
       & 0 & 1
\end{tabular}~.
\end{equa}
\paragraph{Fourth order Runge-Kutta method}
The Kepler problem using forth order Runge-Kutta 
method is basically the same as second order Runge-Kutta, but with 
three  points instead of one point in between to
solve the position and velocity. The tableau we use for this method is \citep{Butcher1963} 
\begin{equa}
\begin{tabular}{c|c c c c}
   0 &  & & &  \\
   1/2 & 1/2 & & &  \\
     1/2 & 0& 1/2  & & \\
      1 & 0& 0  & 1& \\
         \hline
        & 1/6& 1/3  & 1/3&1/6
               \end{tabular}~.
\end{equa}
\subsection{Results for various integration schemes}
\begin{figure}
  \begin{center}
    \includegraphics[width=0.6\textwidth]{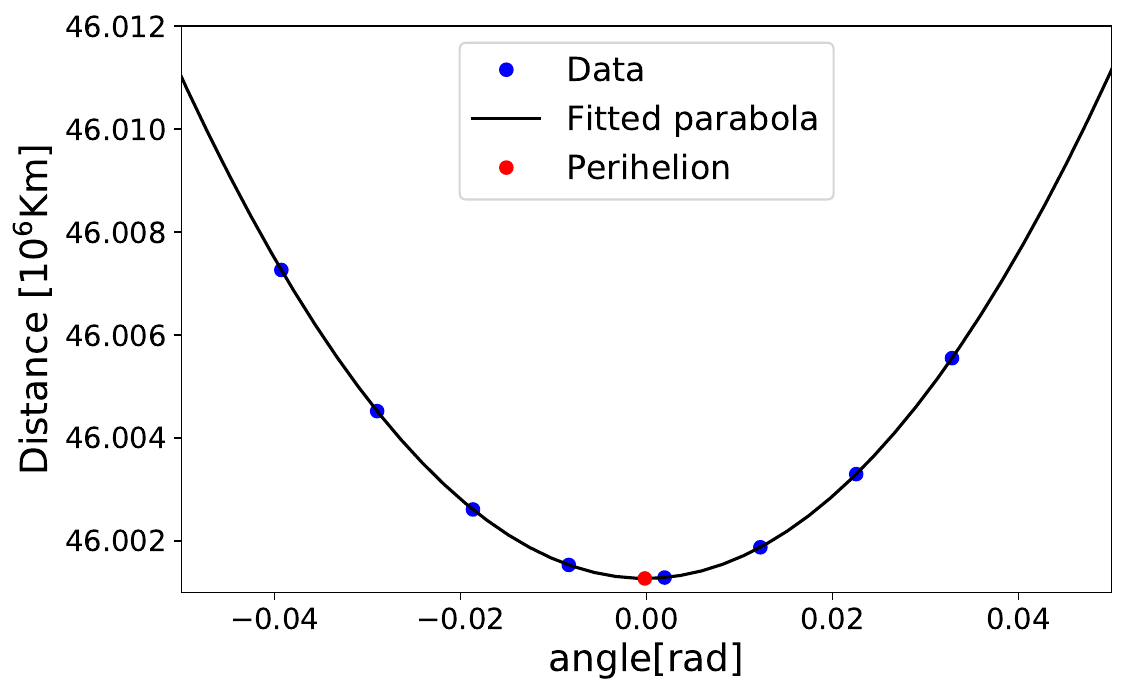}
    \end{center}
  \caption{Detecting the spurious perihelion change for the Kepler
    problem for the Mercury-Sun system. For each time step of integration we determine the angle
    and the distance from the central body (the blue points),
     using the 2nd order Runge-Kutta method.  We then fit a parabola
    through these points, and the minimum of the distance to the Sun
    is the red point (perihelion point). Note that the red point is very slightly to the
    left of 0 and shows the spurious perihelion shift due to the time integration imprecisions.  We use this method to find the perihelion shift and to decrease the errors we take average over three revolutions.}\label{fig:perihelion}
  \end{figure}

\begin{figure}
  \begin{center}
    \includegraphics[width=0.8\textwidth]{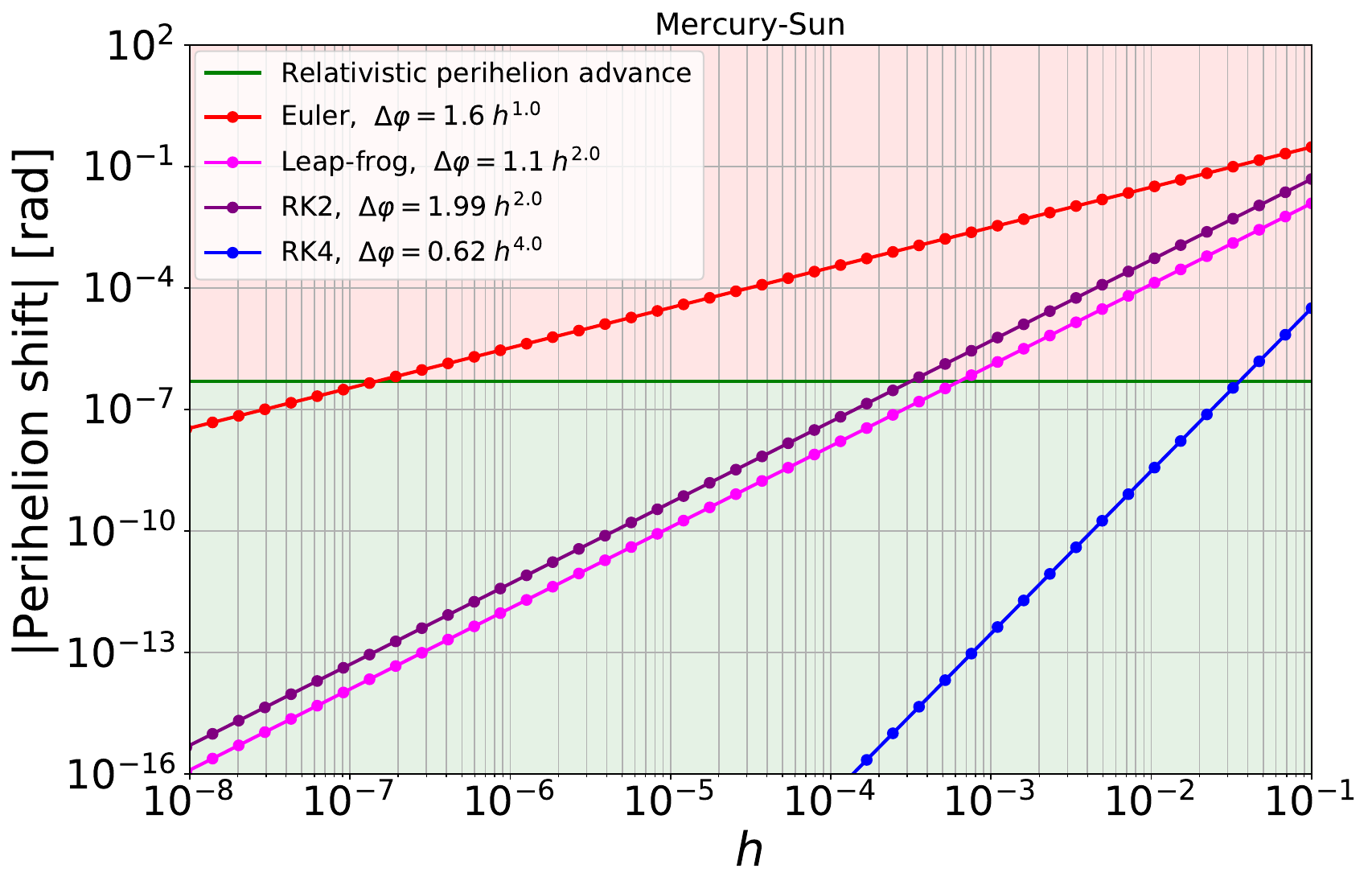}
    \end{center}
  \caption{Achievable precision for different integrators, as a
    function of step size $h$. Shown is the absolute value 
    of the  perihelion shift for the
    Mercury-Sun problem. To make relativistic corrections
    distinguishable, only points in the green region are good
    enough. The data points correspond to $1/3$ of the advance after 3
  rotations.}
\label{fig:timemethod} 
  \end{figure}

 We solve the Kepler problem with the methods described above and find the
 perihelion variation for different time steps $h$. To determine the
 perihelion point of the orbit we choose several points near the minimum distance to the central object after each revolution.
Then we fit a parabola (for Runge-Kutta, we take a 4th order polynomial to
 find the point closest to the central mass) and the minimum of the
 parabola is taken as the perihelion point. \fref{fig:perihelion}
 illustrates how this is done, for the particular example of Mercury.
 The red point is the perihelion.
The spurious shift of the perihelion of Mercury using 2nd order Runge-Kutta method with $h=0.00625$, which is the
case considered in \fref{fig:perihelion}, is $\sim 7.8 \times 10^{-5}$ radians.
We measure the positions in units Giga
meters ($\rm{Gm} \equiv10^9$ m) , time in Mega seconds ($\rm{Ms}=10^6$ sec) and masses in $M_\text{earth} =5.972 \times 10^{24} \, \rm{Kg}$.  In
these units, the initial position (at the perihelion), is 
  $46.001272 (\cos(\theta ),\sin(\theta ))$  where $\theta$ is the
angle between x-axis and semi-major axis. The initial velocity is perpendicular to the line connecting
Mercury and the Sun, with
  magnitude $58.98 $. The potential is $-GM/r$ with $GM= 132733$ measured in the code's units \footnote{In  units ${\rm{G m}^3} \cdot  M_\text{earth}^{-1} \cdot \rm{M}s^{-2}$.}.  When we will
  study the problem on the lattice, the angle $\theta $ will be
  important. In \fref{fig:timemethod} the magnitude of perihelion variation
  for the different time integrators and the step size $h$ is shown. The
  horizontal line shows the value of relativistic perihelion advance,
  the green/red regions respectively show where the time
  integrator precision is/is not good enough to observe relativistic
  perihelion advance.
Because time integrators over- or underestimate the perihelion, we plot
its absolute deviation (which for Newton's law should be zero).
This absolute value sets a limit of how small one has to take a time
step $h$ to be able to detect general relativistic corrections to the orbits.
 
 The relativistic parameter $\Upsilon
 =r_{\rm sch}/ {r_{\rm per}}$ for the Mercury-Sun case with $r_{\rm
   per} =46 \times 10^{6} $ km and $r_{\rm sch} = 2.95 $ km is $\Upsilon
 \approx  6.4 \times 10^{-8}$. The eccentricity of Mercury is  $\approx 0.205630$. Both parameters are considered small as the distance of Mercury to
 the Sun is much larger than the Schwarzschild radius of the Sun.
 We therefore consider also a more extreme case, where the
 relativistic effects are larger, such as the
  stars in the Galactic center known as the S-stars \citep{Parsa} which are
  revolving around the central super massive black hole.
  For one of them, S2, the relativistic parameter at perihelion point is estimated from
  measurements to be
  $\Upsilon \approx 8.8 \times 10^{-5} $, and the eccentricity is
  $0.884$. The details about the relativistic
 parameter and eccentricity can be found in
 Appendix~\ref{parametrization}.

\begin{figure}
  \begin{center}
    \includegraphics[width=0.45\textwidth]{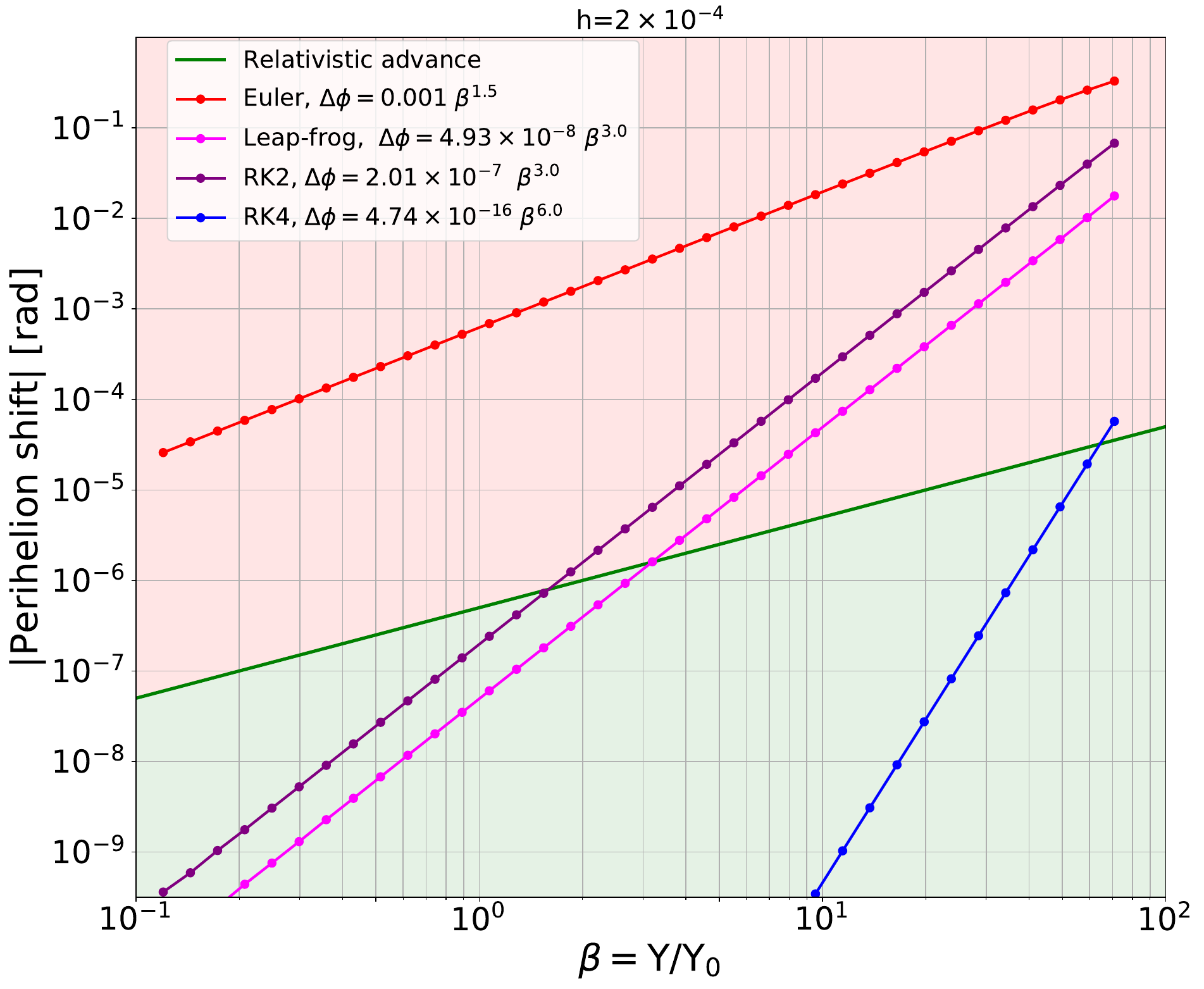}
        \includegraphics[width=0.45\textwidth]{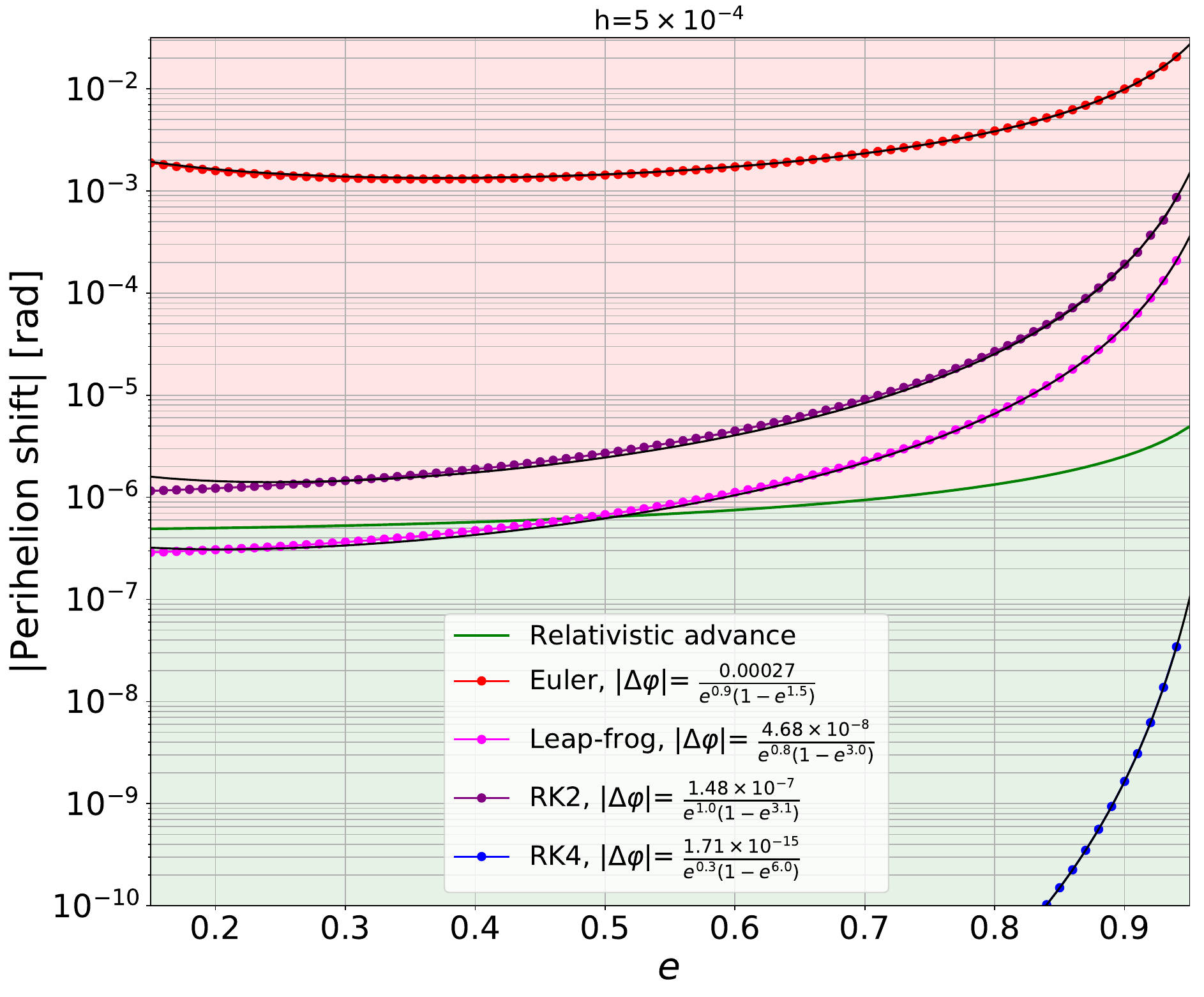}
    \end{center}
  \caption{{\bf Left}: The absolute value of the perihelion shift as a function of the normalized
    relativistic parameter $\beta= \Upsilon/\Upsilon_0$, where $\Upsilon_0$ is the relativistic parameter for Mercury-Sun.  The red region
    shows where the method will fail to discriminate the relativistic
    perihelion advance from the integration errors (for the chosen
    step size of $h=0.0002$). In the green region one can safely
    use the method, for that specific orbit.  When increasing $\beta $, the numerical perihelion shift
    increases, as according to \eqref{rel_ic} the velocity of the object in the perihelion point scales like $\sqrt{\Upsilon}$, while the perihelion distance scales like $\frac{1}{\Upsilon}$. In all the methods the slope of the curve is higher
    than the slope of the relativistic advance curve $\Delta \varphi
    \sim \Upsilon$, which shows that for the orbits with large
    relativistic parameters, one has to choose the method and the time
    step very carefully. {\bf Right}: The same representation
    as a function of eccentricity $e$. In all the methods, by increasing the eccentricity the
    numerical perihelion variation increases, as according to the \eqref{ecc_ic} the velocity of the object and perihelion distance rescales respectively by $ \sqrt { \frac{1+e}{1+e_0}}$ and
$ {\frac{1-e}{1-e_0}} $. In order to be able to
    measure the relativistic perihelion advance at each eccentricity
    we need to use the method with the appropriate step size, for
    example Euler and second order Runge-Kutta do not work for any
    eccentricity, while leapfrog is good for $e\lesssim 0.5$ and
    fourth order Runge-Kutta works perfectly for all
    eccentricities. All data points correspond to $1/3$ of the advance after 3
  rotations
 }\label{fig:advance_upsilon} \label{fig:all_methods-eccentricity} 
  \end{figure}

 In our numerical study, we cover therefore a large range of these
 parameters. Our results are summarized for the four integration
 methods in \fref{fig:timemethod} as a function of the time step
 $h$. In \fref{fig:advance_upsilon} the comparison is done as a
 function of
 $\Upsilon/\Upsilon_0$, where $\Upsilon_0$ is the relativistic
 parameter at perihelion point for the Mercury-Sun system. We also
 show the dependence on eccentricity. 

  \section{Force interpolation} \label{discretization}

  \begin{figure}
  \begin{center}
    \includegraphics[width=0.4\textwidth]{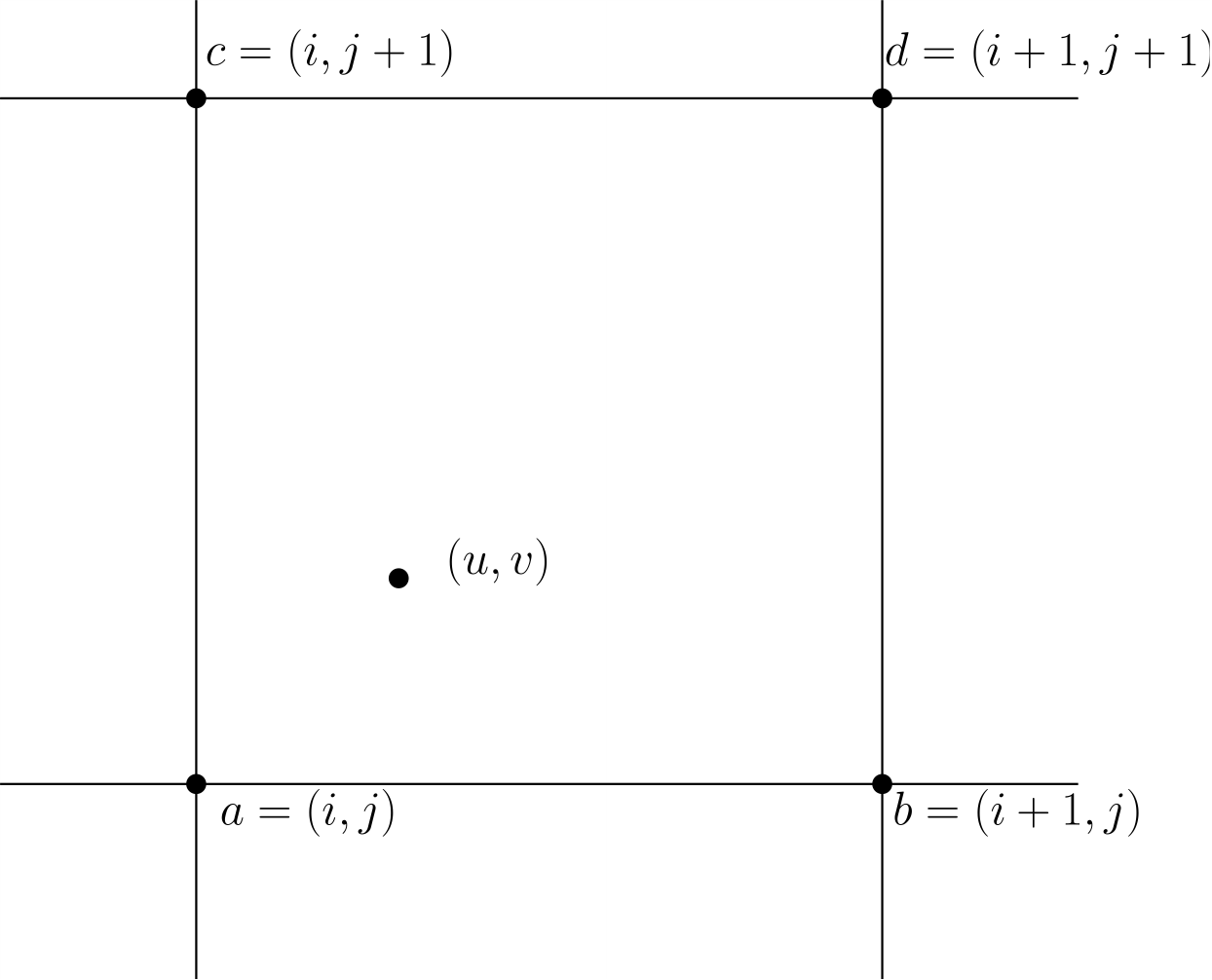}
    \end{center}
  \caption{Notation for a plaquette on the unit mesh in $\real^2$.}\label{fig:1}
  \end{figure}

 Having considered the numerics of the classical methods, we now study
 the effect of discretizing space. We again restrict attention to two
 dimensions and set, throughout, the lattice spacing equal to
 $dx$.\footnote{Due to the discretization, the angular momentum
   vector might not be conserved and we might have 3D motion, here we assume that the force perpendicular to the plane of motion vanishes.}
In particular, we study  the two force interpolations (linear and
bilinear) which are mainly used in N-body simulations, see {\it e.g.},
\citep{Springel2005} and \citep{Adamek2016}, and for which we will
present numerical results.
A very useful systematic derivation of finite elements for derivatives
and differential complexes can be found in \citep{arnold2002}. The setup is as follows:
We are given a potential $\Phi$, in our case the Newtonian potential $\Phi(x_1,x_2)=-GM/\sqrt{x_1^2+x_2^2}=-GM/r$, from which we want to derive the
forces on the particles.
In the bilinear (quadratic) method the mesh is given by integer coordinates (in
$\integer^2$),  and we assume that $\Phi$ is known in all points
$(i,j)$, with $i,j\in \integer$.\footnote{Finite elements are of course obtained more easily
  on triangular lattices, but, because of requirements of large
  parallel computations, we study the lattice $\integer^2$.}
The force at lattice point $i,j$ is then approximated by a vector with components
\begin{equa}
  f^{(x)}_{i,j} &= \frac{\Phi_{i+1,j}-\Phi_{i-1,j}}{2}~,\\
  f^{(y)}_{i,j} &= \frac{\Phi_{i,j+1}-\Phi_{i,j-1}}{2}~.\\
\end{equa}
Note that the difference is taken over 2 mesh points around the point
of interest.
Assume that the point $x=(u,v)\in\real^2$ lies in the square with
corners
\begin{equa} 
 a=(i,j)~,\quad b=(i+1,j)~,\quad c=(i,j+1)~,\quad d=(i+1,j+1)~,
\end{equa}
cf.~\fref{fig:1}.
We let $f_a^{(x)}=f_{i,j}^{(x)}$, and similarly for the other corners
and the direction $(y)$.

Let $i$ be the integer part of $u$ and let $j$ be the integer part of
$v$, and set $\xi=u-i$, $\eta=v-j$.
The interpolated forces are then given by
\begin{equa}[eq:springel]
  F^{(x)}(u,v)=&\left(f^{(x)}_a (1-\xi)+f_b^{(x)}\xi \right)\cdot
  (1-\eta)
  +\left(f^{(x)}_c (1-\xi)+f_d^{(x)}\xi \right)\cdot \eta~,\\
  F^{(y)}(u,v)=&\left(f^{(y)}_a (1-\eta)+f_c^{(y)}\eta \right)\cdot
  (1-\xi)
  +\left(f^{(y)}_b (1-\eta)+f_d^{(y)}\eta \right)\cdot \xi~.\\
\end{equa}
Note that ``$c$'' and ``$b$'' change position between the $x$ and $y$
components.This interpolation method is called {\bf{bilinear method}} as it is combination of two linear interpolations along the square, so it is a quadratic interpolation \citep{arnold2002}.
 This interpolation is
continuous across the boundaries in both directions and for both
components of the vector field. 
To verify this, one can for example restrict to the line connecting
$a$ and $b$. Then $\eta=0$, and therefore one gets
\begin{equa}
  F^{(x)}(u,v)=&f^{(x)}_a (1-\xi)+f_b^{(x)}\xi~,\\
  F^{(y)}(u,v)=&f^{(y)}_a 
  (1-\xi)+f^{(y)}_b  \xi~.\\
\end{equa}
The important thing is that the values only depend on $a$ and $b$, but
not on $c$ and $d$ and so continuity is guaranteed. The 3 other edges
are similar.

In this scheme, as is well known, one needs 8
evaluations of $\Phi$ per plaquette. When the mesh size is $dx$ instead of $1$, all
the calculations scale accordingly.

The other method, {\bf{linear method}} which is also widely used in
cosmological N-body simulations {\it e.g.}, \citep{Adamek2016} is given, with 
similar notation---using $g$ and $G$ instead of $f$ and $F$---by 
\begin{equa}
  g^{(x)}_{i,j}&=\Phi_{i+1,j}-\Phi_{i,j}~,\\
  g^{(y)}_{i,j}&=\Phi_{i,j+1}-\Phi_{i,j}~,\\
  G^{(x)}(u,v)&=g^{(x)}_a(1-\eta)+g^{(x)}_c\eta~,\\
  G^{(y)}(u,v)&=g^{(y)}_b(1-\xi)+g^{(y)}_d \xi~,
\end{equa}
with $a,\dots,d$ as before.
This method is of lower order than the previous one, and needs fewer
evaluations. The advantage is that they need less memory, but of
course, it is only $1^{\rm st}$ order.

Note that if $(u,v)$ crosses the line connecting $a$ and $b$, then
$G^{(x)}$ is continuous, but $G^{(y)}$ has a jump discontinuity  (of
order about $\OO(\Phi_{i,j}-\Phi_{i+1,j})$ when $i$ and $j$ are not too close to 0).
Similar considerations hold on the other boundaries of the unit
plaquette.

This scheme only needs 4 evaluations of $\Phi$ per 2-dimensional plaquette, but the
interpolation is not continuous. The Kepler problem can still be
integrated numerically, but there will appear a spurious phase
shift which is caused by the discontinuity. But the
numerical errors again scale with the mesh size, albeit on a larger
scale than in the first method.

We will now present the numerical results for these cases, and then
discuss the limitations they imply on trajectories in N-body simulations. Of course, often calculations are done in
$\real^3$, resp.~$\integer^3$, but for the study of numerical issues,
2 dimensions are enough. Restriction to 1 dimension is too easy, since
the two methods coincide in that case.

\section{Discretization vs relativistic perihelion advance}

We have seen that high precision is needed to discriminate
relativistic effects in the planar two-body problem. As several
N-body codes use---in addition to the standard numerical integration
schemes, a \emph{discretization} of space---we now study the effects of
these discretizations. To concentrate on them, we use a numerical
integration of very high precision (ODEX, tolerance $8\cdot10^{-11})$ so that
the effects described earlier are minimal, and the effect of
discretization becomes visible.


As we want to measure the perihelion advance due to the discretization we stick to the general equations in which we do not use the symmetries of the Kepler problem. We just assume that the motion is on a plane and we solve the equations for the relative distance between two masses assuming that $m \ll M$,
  \begin{equa}
(\dot x, \dot y) &= (v_x, v_ y)~,\\
(\dot{v}_x, \dot{v}_ y)&= ( H^{(x)}(u,v) , H^{(y)}(u,v))~,
  \end{equa}
where $ H^{(x)}(u,v)$ is $F^{(x)}(u,v)$ or $G^{(x)}(u,v)$ as defined in Section~\ref{discretization} .

We have, so far, analyzed in detail how much numerical precision is
needed to detect general relativistic effects, as a function of eccentricity and the relativistic parameter $\Upsilon$.

\begin{figure}[ht!]
\begin{center}
    \includegraphics[width=0.45\textwidth]{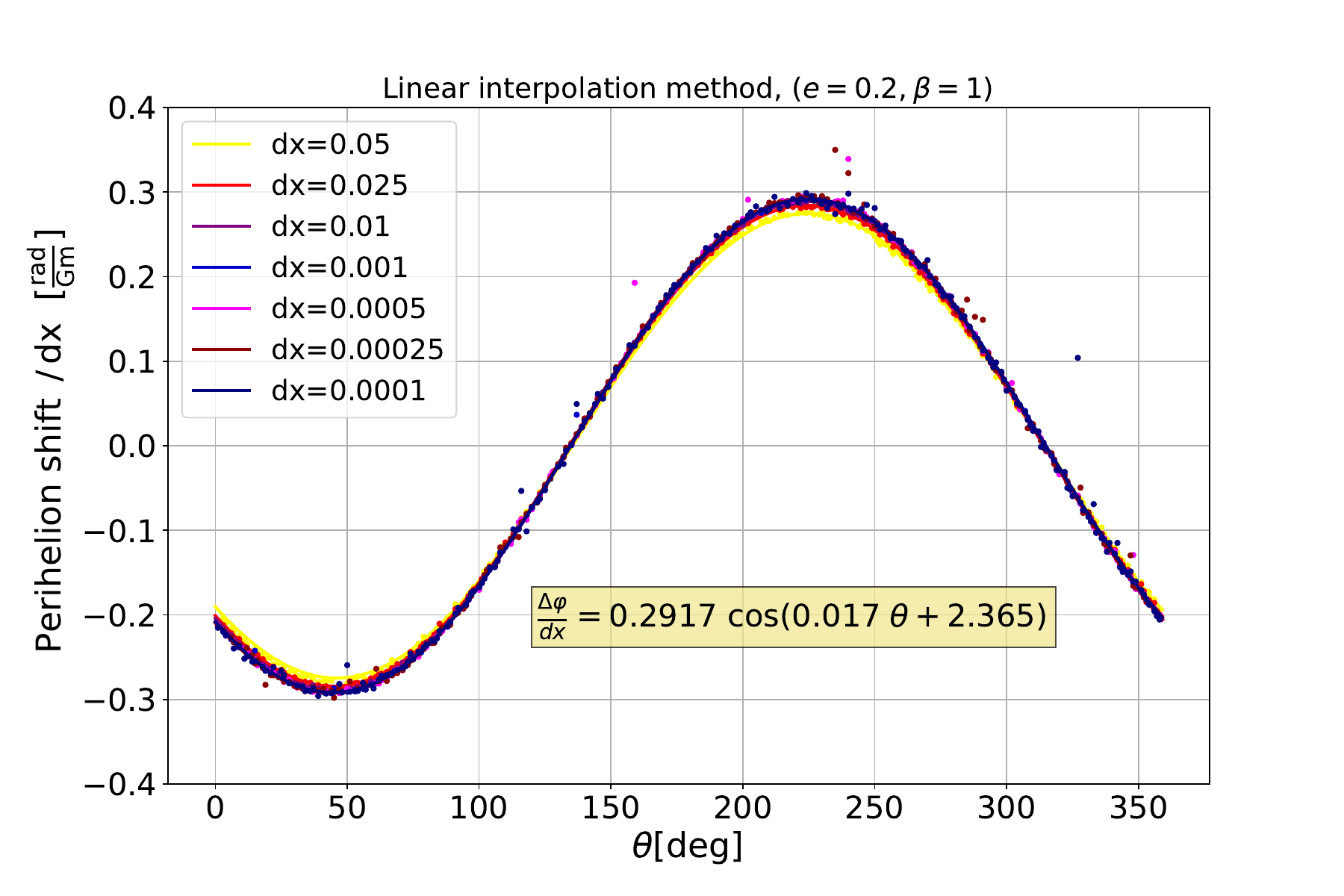}
   \includegraphics[width=0.45\textwidth]{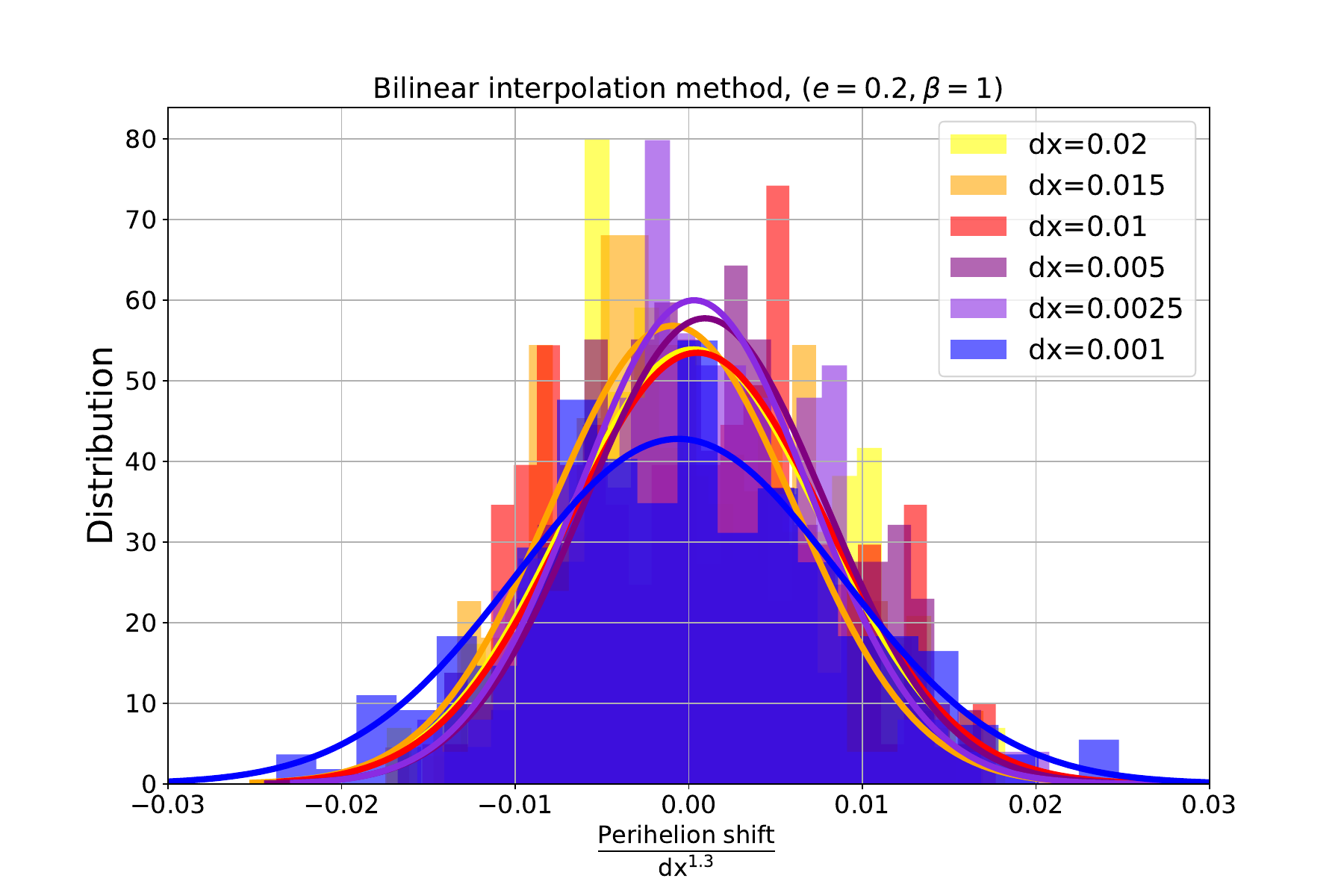}   \end{center}
  \caption{{\bf Left}: Linear interpolation: The advance of the
    perihelion due to the discretization effect depends  on the initial angle of the perihelion. For each
    value of $dx$, we show the deviation in radians,
  divided by $dx$. The curves clearly coincide. This shows that the
  deviations scale linearly with $dx$.\\
  {\bf Right}: Bilinear interpolation: 
We consider the perihelion shifts, divided by $dx^{1.3}$, for 180
equally spaced initial angles of the orbit. The bar graphs show the
distribution of these quantities, for various choices of $dx$. We see
that they obey a Gaussian fit (the solid lines). This shows that the shifts
are random.}\label{fig:perihelionlinear}
  \end{figure}

\begin{figure}[ht!]
\begin{center}
    \includegraphics[width=0.7\textwidth]{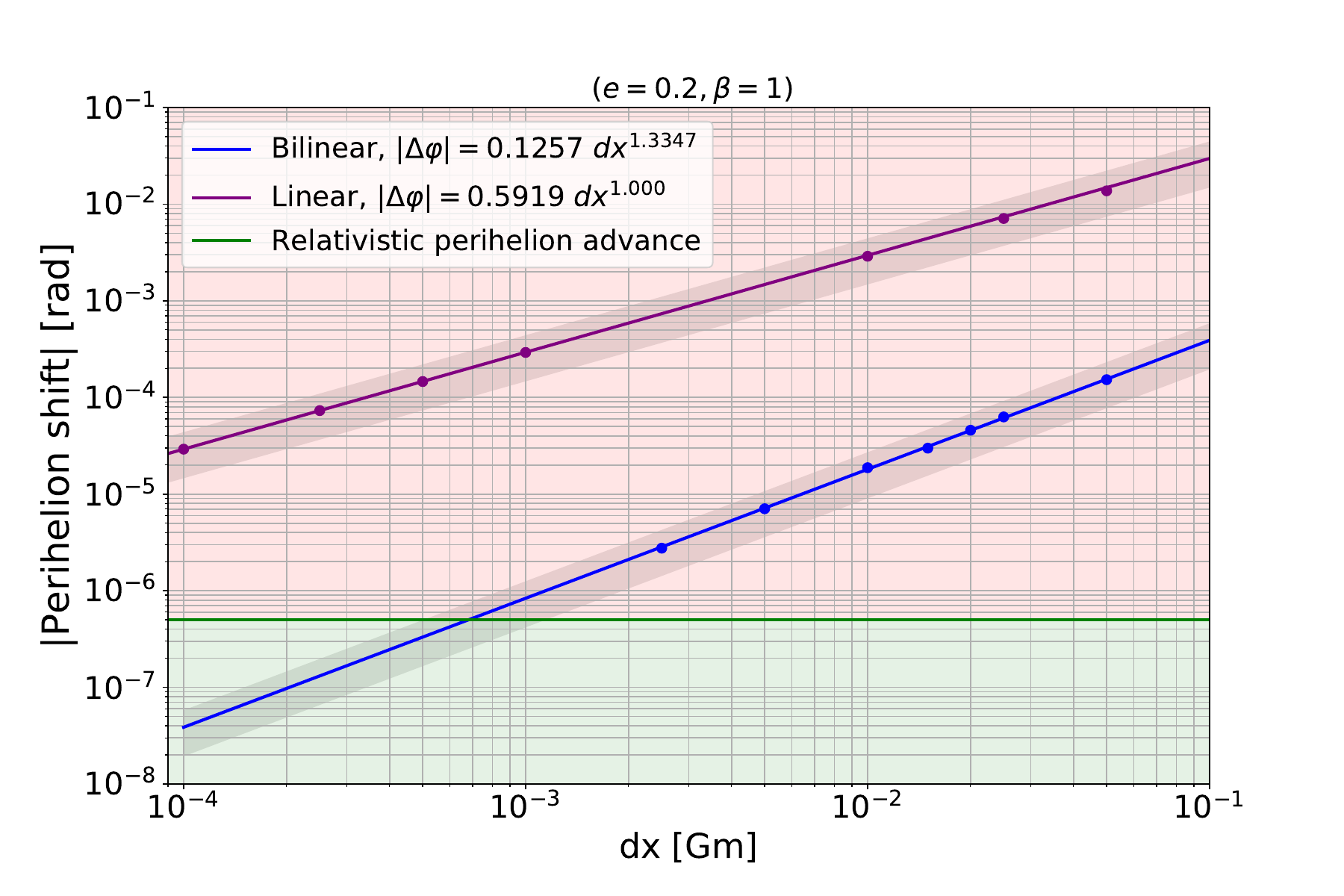}
  \caption{Dependence of the standard deviation of the
    perihelion shift for the two methods, as a function of
    $dx$. For the bilinear method, numerical fluctuations are too
    large to get reliable results for $dx\lesssim 10^{-3}$. The grey regions around the fitted blue and magenta lines show the $50\%$ deviation from the central fit value.}\label{dxmercury}
  \end{center}
\end{figure}

We start by presenting results for the bilinear interpolation
\eref{eq:springel}. The simulations are done as follows: We take the
parameters for Mercury, with initial position at $46.001\cdot
(\cos(\theta ),\sin(\theta ))$ and a velocity perpendicular to the
Mercury-Sun line, of
magnitude $58.98$, in the counterclockwise direction. 
We require a tolerance of $8 \times 10^{-11}$, which is attainable with
quadruple precision, using ODEX. For each value of $\theta$, we
determine the time for 1, 2, and 3 returns to the perihelion. 
The
perihelion is found by looking for that angle where the distance from
the sun is minimal. This angle is
found by linear and quadratic bisection, up to machine precision,
using the ``continuous output'' from ODEX.\footnote{We use the standard
  algorithm ``zeroin'' of Dekker.} We repeat this  for 180
initial angles covering 360 degrees in steps of 2 degrees, and this gives us $3\times 180$
data points.\footnote{For example, after 3 turns, we divide the total
  angle by 3, we do not take the difference between the angle for 3
  and 2 turns. Of course, errors on these points will average out
  somewhat.}
In \fref{fig:perihelionlinear} and \fref{dxmercury} we show the results for several values of
$dx$, for linear and bilinear approximation.
Further inspection shows that these distributions are close to
Gaussian, but the variance is somewhat smaller than $dx$, actually
$dx^{1.3}$ is a reasonable approximation. As we mentioned before, an
analytic estimate of this variance is difficult, because ODEX works
with variable step size and order, with quite dramatic changes near
the edges of the plaquettes.

\begin{figure}[ht!]
\begin{center}
    \includegraphics[width=0.45\textwidth]{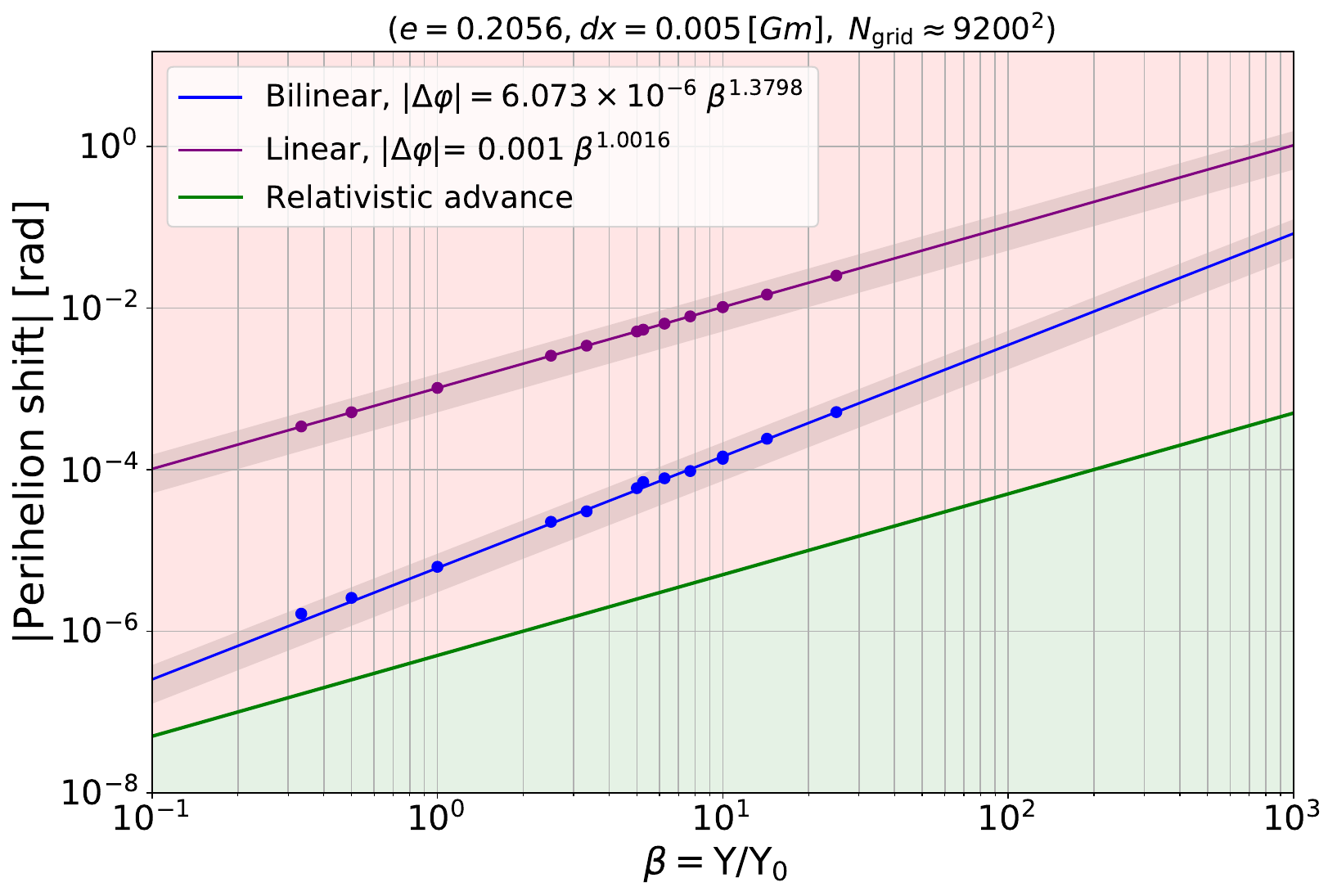}
   \includegraphics[width=0.45\textwidth]{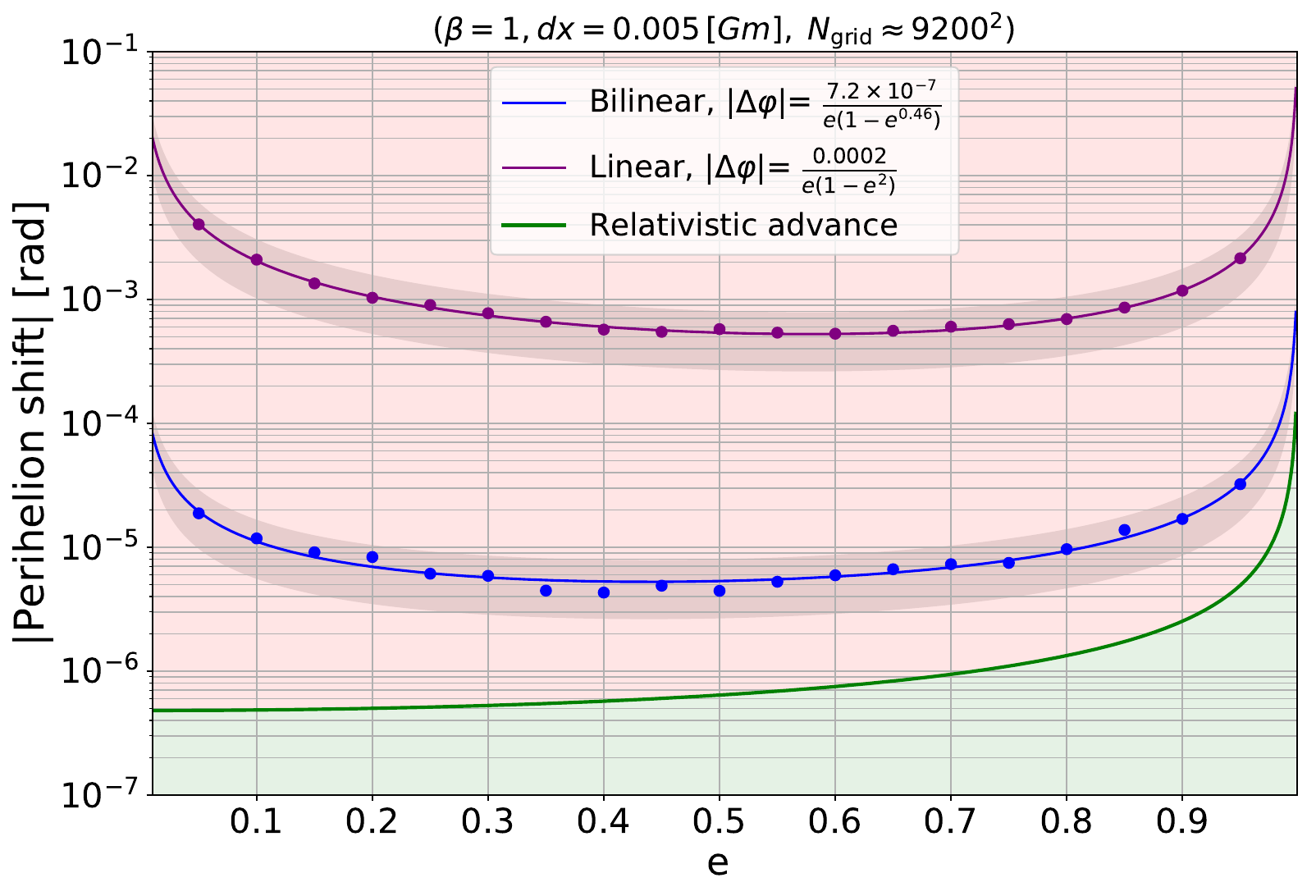}   \end{center}
  \caption{Comparison of linear vs.~bilinear interpolation. {\bf
      Left}: Behavior as a function of the theoretical relativistic parameter
    $\beta $. The green line shows the relativistic perihelion advance
    and therefore the green/magenta area determines the regions for
    $\beta $ where one can/cannot detect the perihelion advance. Note
  that for the parameters $e$, $dx$, and $N_{\rm grid}$ used in the
  figure, no method is able to detect relativistic
  corrections.\\
  {\bf Right:} The same study, now as a function of eccentricity
  $e$. According to the green/magenta regions the relativistic effects
  can not be detected for the choice of $\beta $, $dx$, and  $N_{\rm
    grid}$. 
  }\label{fig:betaecc}
  \end{figure}

In the case of the linear interpolation, the
discontinuity leads to an effective advance of the perihelion, which
furthermore depends strongly
on the initial angle $\theta $. Qualitatively, this can be understood
by the angles at which the orbit crosses the discontinuities.
Using otherwise the same
parameters as above, the results are summarized in \fref{fig:betaecc}.
The advance $A$ of the perihelion follows closely a cosine (with a
phase-shift) $A/dx\sim 0.145\cos(\theta+2.34)\sim 0.145 \cos(\theta+3\pi/4)$. We also
checked that the advance of the perihelion changes sign if the initial
velocity changes sign. Also note that the average of the advance of the
perihelion is close to zero.
To generalize the results to include different orbits in N-body
simulations especially the ones with high eccentricity and high
relativistic parameter, we studied how the perihelion shift depends on
eccentricity and relativistic parameter.
These are shown in \fref{fig:betaecc}. Given the number of particles
considered in current N-body simulations ({\it e.g.}, $7000^3$ in
\citep{Yu2017} ) and the restrictions of current hardware, (same
number of lattice points), we see that relativistic corrections of the
orbits can not be detected.

An interesting effect of the discretization for fixed $dx$ is the
dependence on eccentricity $e$. We observed that the deviations scale
about as $\frac{1}{e(1-e^2)}$ for the linear interpolation. This means
that the effect is largest at extreme values of $e$. The deviation is
also proportional to $\Upsilon$, while the relativistic correction is
proportional to  $\frac{\Upsilon}{1-e^2}$.

The code for such tests and for other parameters can be
obtained from the authors.

\section{Conclusions}

We see that the numerical study of relativistic effects can have two
problems. First, the integration method must choose a small enough time
step to reach a precision which is better than the size of the relativistic
correction. Second, if, additionally, the forces are discretized, the
grid size must be quite fine, so that the relativistic corrections are
not washed out by the approximation.

In particular, our results allow one to estimate for which choices of
$\Upsilon$, $e$, and $\theta$, the relativistic effects are larger than
the numerical and discretization effects generated by $h$ and $dx$.

\section*{Acknowledgements}
We thank Martin Kunz for asking questions leading to this paper, and
for helpful discussions. We thank Ruth Durrer and Jacques Rougemont
for helpful comments about our manuscript. JPE acknowledges partial support
by an ERC advanced grant ``Bridges'' and FH acknowledges financial support from the Swiss National Science Foundation.

\appendix 
\section{Newtonian and relativistic orbits}
Here first we review the equations of motion and properties of the orbits in the mechanics of Newtonian particles, then we derive the equations for general relativistic motion.
In 2-dimension polar coordinates $(r,\phi)$ the Newton equations take
the form
\begin{equa}
\ddot{r} - r \dot{\phi}^2 =  -\frac{G M}{r^2}
\end{equa}
Writing the angular momentum per mass $l$ in polar coordinates results in
\begin{equa}
l = r^2 \dot{\phi}~.
\end{equa}
Changing the variable to $u (\phi)=  1/ {r (\phi)}$ gives
\begin{equa}
\frac{d^2}{d\phi^2}u  + u = \frac{G M }{l^2}~. \label{newt_ode}
\end{equa}
To obtain the relativistic perihelion advance we repeat, for the
convenience of the reader, some parts of
 \citep[p.~193]{Stephani:1982ac}. The Schwarzschild metric as a spherically symmetric vacuum solution reads,
\begin{equa} 
ds^2=-\Big(1-  \frac{r_{\rm sch} }  {r} \ \Big) c^2dt^2  +  
\Big( {1-   \frac{r_{\rm sch}}  {r} }  \Big)^{-1}{dr^2}   + r^2 \big( d \theta^2 + \sin ^2 \theta \, d \phi ^2 \big)~,
\end{equa}
where $\Rsch$ is the Schwarzschild radius, defined by  $r_{\rm sch}=
\frac{2 GM}{c^2}$, $t$ and $r, \theta, \phi$ are respectively time and
spatial spherical coordinates.\footnote{In this section, $\theta$ is
  not the angle of the major axis, but just one of the 3 Euler coordinates.} To obtain geodesic equations one starts
from the classical action of massive test particle, 
\begin{equa}
\mathcal{A}=-m_o  \int \sqrt{-g_{\mu \nu} \dfrac{\d x^{\mu}}{d \tau}
  \dfrac{\d x^{\nu}}{d \tau} } d \tau~,
\end{equa}
where $m_o$ is the mass of the object. Applying Euler-Lagrange
equation to the Lagrangian of the test particle gives four equations,
in which one sees that  the angular momentum is conserved. Therefore,
the motion is in a plane. By
simple algebra on the equations one finds the equation of motion as
\begin{equa}
\frac{d^2}{d\phi^2}u + u = \frac{G M }{l^2}+  \frac{3}2 r_{\rm sch} u^2 \label{gr-ode}~,
\end{equa}
which is like Newton's equation \eref{newt_ode} plus a term which is
coming from relativistic correction. We solve \eref{gr-ode} to obtain
the relativistic perihelion advance per period, which is well approximated by  
\begin{equa}\label{eq:advance}
\Delta \phi_p \approx \frac{3 \pi r_{\rm sch}}{a \big(1- e ^2 \big)}~,
\end{equa}
where $a$ is semi-major axis and $e$ is the eccentricity of the
orbit.

For Mercury, this leads to the well-known advance of 42.98 arc sec
perihelion advance per century, or $\sim 0.103$ arc sec per period.
\label{sec:maths} 

\section{The parameterization of orbits} \label{parametrization}
To see the effect of discretization on different orbits in N-body simulations,
we parameterize a general orbit with three parameters $(\Upsilon, \theta, e)$, where $e$ is
the eccentricity, $\Upsilon$ is the relativistic parameter at perihelion and $\theta$ is the angle
of semi-major axis with the lattice squares. It is important to note that, these parameters are enough to explain any closed orbits in N-body simulations. Moreover having the three parameters one could uniquely construct the mass of central object as well as the initial position and velocity of the particle.

\begin{figure}[h!]
  \begin{center}
    \includegraphics[width=0.6\columnwidth]{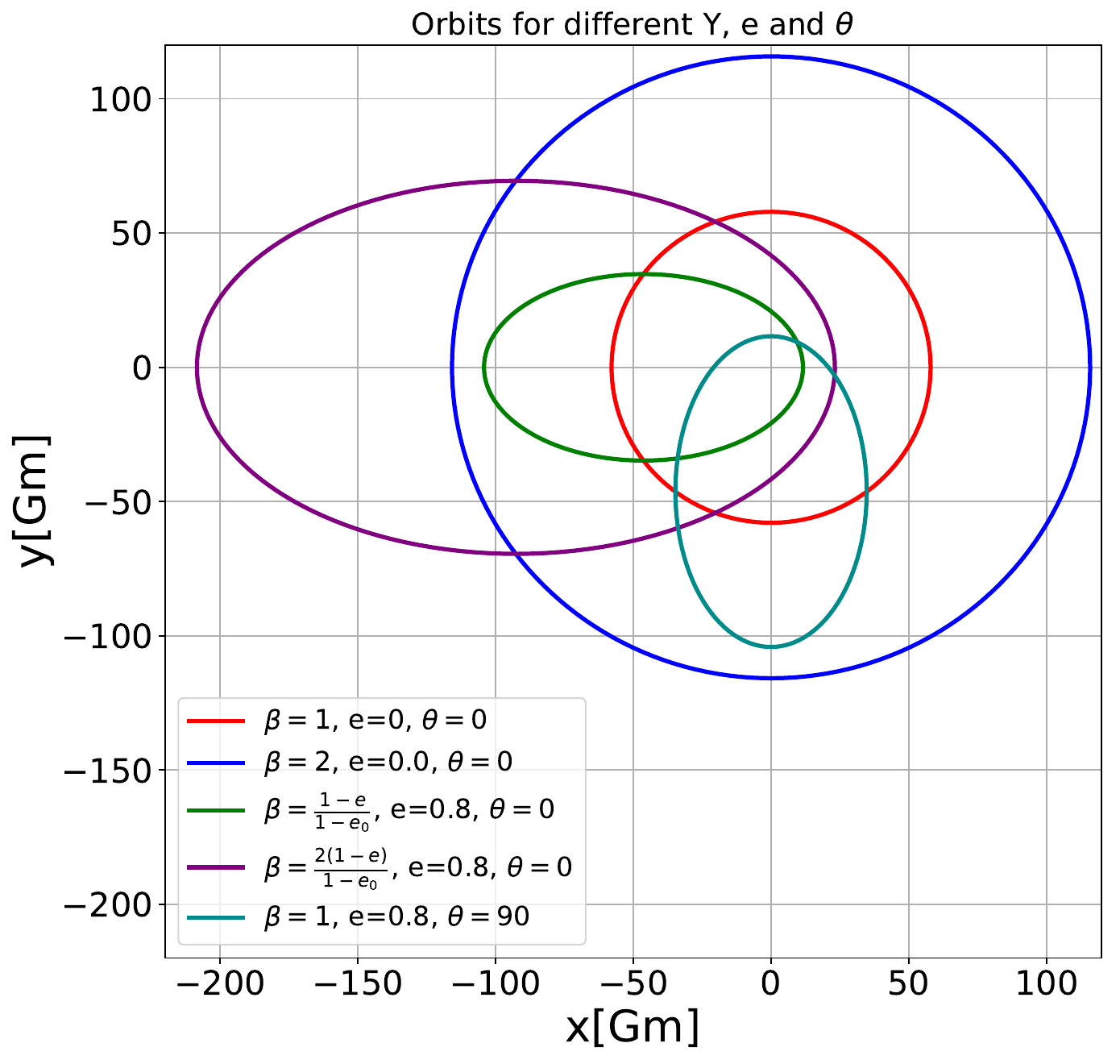}
    \end{center}
     \caption{Some examples of the parameterization of elliptic
       orbits, which show the role of $\beta = \Upsilon/\Upsilon_0$, e, and $\theta$. The orbits are obtained by solving the differential equations. }     \label{fig:orbits}
 \end{figure}

\paragraph{Relativistic parameter}
The relativistic parameter $\Upsilon = \frac{r_{\rm sch}}{r_{\rm per}}$, for a fixed mass of central object it shows the scale of the orbits and for a fixed  size of the orbit it is an indicator of the mass of central object. If we assume that the mass of central object is fixed, by changing the relativistic parameter, different quantities of the orbit would scale as following,
\begin{equa}[3] \label{rel_ic}
M &\rightarrow  M~,\qquad &r_{\rm per} &\rightarrow   \frac{r_{\rm per}}{\Upsilon}~, \\
 T &\rightarrow \Upsilon^{3/2} {T}~,\qquad 
 &v_{\rm per} &\rightarrow {\sqrt{\Upsilon}} {  v_{\rm per}  } ~.
\end{equa}
$r_{\rm per}$ is the  perihelion radius, T is the period of the orbit and $ v_{\rm per} $ is the velocity of the object in the perihelion point. To rescale the orbit for the fixed central body mass and fixed eccentricity one has to change the initial conditions as following to obtain the new orbit,
\begin{equa}[3 ]
x_0&=&  \frac{r_{\rm per}}{\Upsilon}~,\qquad
y_0&=&0~, \\
v_x&=& 0~,\qquad
v_y&=&{\sqrt{\Upsilon}} {v_{\rm per}}~.
\end{equa}
We could of course change the central object mass instead of changing the size of the orbit while having the same relativistic parameter. 
\paragraph{Eccentricity}
Another parameter which is important in characterizing an orbit is the eccentricity, to change the eccentricity we keep the
semi-major axis length  fixed and we change the positions and
velocities in the perihelion point to recover the desired eccentricity for the orbits
\begin{equa} \label{ecc_ic}
  r_{\rm{per}} &\to r_{\rm{per}} \, {\frac{1-e}{1-e_0}} ~,\qquad
v_{\rm{per}} &\to   v_{\rm{per}} \sqrt { \frac{1+e}{1+e_0}} ~,
\end{equa}
where $e$ is the new eccentricity and  $e_0$ is the reference eccentricity (in our case mercury). Note that changing eccentricity also results in changing the perihelion distance and relativistic parameter.
\paragraph{Rotation}
It appears that the angle between the semi-major axis and the lattice squares, is an important parameter specially in the linear force interpolation.  To rotate the orbit by  angle $\theta$ we can follow the coordinate transformations and start from the following initial condition to obtain the correct orbit,
\begin{equa}
x_0 &=&  r_{\rm{per}} \cos (\theta)~,\qquad
y_0 &=&  r_{\rm{per}}  \sin (\theta) ~, \\
v_x &=&    -v_{\rm{per}}   \sin (\theta)~,\qquad
v_y &=& v_{\rm{per}} \cos (\theta)~.
\end{equa}
In \fref{fig:orbits} we have illustrated the orbits with different ellipticity, relativistic parameter and angle obtained from numerical results.

\end{document}